\documentclass{amsart}

\usepackage{amsthm}
\usepackage{amsfonts}
\usepackage{amsmath}
\usepackage{amssymb}
\usepackage{mathrsfs}
\usepackage{fullpage}
\usepackage{stmaryrd}
\usepackage{hyperref}
\usepackage{verbatim}

\theoremstyle{plain}

\theoremstyle{definition}

\theoremstyle{remark}

\newcommand{\Begin}[2]{\begin{#1}\label{#2}}

\newcommand{\bbQ}{\mathbb{Q}}
\newcommand{\bbR}{\mathbb{R}}

\newcommand{\CA}{\mathcal{A}}
\newcommand{\CB}{\mathcal{B}}
\newcommand{\CM}{\mathcal{M}}
\newcommand{\CN}{\mathcal{N}}

\newcommand{\SCRL}{\mathscr{L}}
\newcommand{\SCRU}{\mathscr{U}}

\newcommand{\lborel}{{\Delta_1^1}}
\newcommand{\cantorspace}{{{}^\omega 2}}
\newcommand{\bairespace}{{{}^\omega\omega}}

\newcommand{\KP}{\mathsf{KP}}

\newcommand{\SR}{\mathrm{SR}}
\newcommand{\CC}{\mathcal{C}}
\newcommand{\RR}{\mathrm{R}}

\newcommand{\EF}{\mathrm{EF}}

\begin{document}

\title{Bounds on Scott Ranks of Some Polish Metric Spaces}

\author{William Chan}
\address{Department of Mathematics, University of North Texas, Denton, TX 76203}
\email{William.Chan@unt.edu}

\begin{abstract}
If $\CN$ is a proper Polish metric space and $\CM$ is any countable dense submetric space of $\CN$, then the Scott rank of $\CN$ in the natural first order language of metric spaces is countable and in fact at most $\omega_1^\CM + 1$, where $\omega_1^\CM$ is the Church-Kleene ordinal of $\CM$ (construed as a subset of $\omega$) which is the least ordinal with no presentation on $\omega$ computable from $\CM$.

If $\CN$ is a rigid Polish metric space and $\CM$ is any countable dense submetric space, then the Scott rank of $\CN$ is countable and in fact less than $\omega_1^\CM$.
\end{abstract}

\thanks{June 10, 2019. The author was supported by NSF grant DMS-1703708.}

\maketitle

%\tableofcontents

\section{Introduction}\label{introduction}

A common task in mathematics is to distinguish different mathematical structures subjected to the restriction of various first order languages. The Scott analysis is a general model theoretic concept that attempts to find an $\SCRL$-isomorphism invariant of an $\SCRL$-structure $\CM$, where $\SCRL$ is a first order language. Informally, if two tuples $\bar{a}$ and $\bar{b}$ of $\CM$ of the same length can be distinguished from each other by an infinitary $\SCRL$-formula, the Scott analysis would attempt to assign an ordinal that indicates how difficult it is to distinguish these tuples.

The Scott rank of tuples can be defined by the back-and-forth relations (see Definition \ref{scott rank}): Let $\bar{a} = (a_0,...,a_{p - 1})$ and $\bar{b} = (b_0,...,b_{p - 1})$ be tuples of length $p$ from an $\SCRL$-structure $\CM$. 

One says that $\bar{a} \sim_0 \bar{b}$ if and only if the map taking $a_i$ to $b_i$ for  $i < p$ is a partial $\SCRL$-isomorphism of $\CM$ into $\CM$.

Assume $\sim_\alpha$ has been defined, one says that $\bar{a} \sim_{\alpha + 1} \bar{b}$ if and only if for all $c \in \CM$, one can always find a $d \in \CM$ so that the elongated tuples satisfy the relation $\bar{a}c \sim_\alpha \bar{b}d$ and similarly in the other direction with the role of $\bar{a}$ and $\bar{b}$ reversed.

Assume that $\sim_\beta$ has been defined for all $\beta < \alpha$, then one defines $\bar{a} \sim_\alpha \bar{b}$ if and only if for all $\beta < \alpha$, $\bar{a} \sim_\beta \bar{b}$. 

If there is an $\alpha$ so that $\neg(\bar{a} \sim_\alpha \bar{b})$, then using the wellfoundedness of the class of ordinals, $\SR(\bar{a},\bar{b})$ is defined to be the minimal such ordinal. Otherwise, one will say $\SR(\bar{a},\bar{b}) = \infty$. Intuitively, $\SR(\bar{a},\bar{b}) = \infty$ indicates that the two tuples are indistinguishable by an infinitary $\SCRL$-formula. If $\SR(\bar{a},\bar{b}) \neq \infty$, then $\SR(\bar{a},\bar{b})$ is an ordinal measuring how difficult it is to distinguish these two tuples. For instance, $\SR(\bar{a},\bar{b}) = 0$ means $\neg(\bar{a} \sim_0 \bar{b})$. Thus there is an atomic formula that evaluates differently between $\bar{a}$ and $\bar{b}$. $\SR(\bar{a},\bar{b}) = 1$ would mean that atomic formulas can not distinguish $\bar{a}$ and $\bar{b}$, but there is a formula consisting of an existential quantifier over an atomic formula that evaluates differently between $\bar{a}$ and $\bar{b}$. In this way, the Scott ranks of tuples are closely relate to the ranks of an infinitary $\SCRL$-formulas that can be used to distinguish tuples. By taking supremum of all possible pairs of tuples of the same length (varying over all possible lengths), one obtains an ordinal for the entire structure $\CM$, called the Scott rank of $\CM$.

Another useful perspective on the Scott rank of tuples $\bar{a}$ and $\bar{b}$ is through a two player game called the Ehrenfeucht-Fra\"{i}sse game $\mathrm{EF}_\alpha^{\CM,\bar{a},\bar{b}}$, where $\alpha$ is an ordinal. Player 1 at each turn plays a pair $(\beta,x)$ where $\beta < \alpha$ is less than any previous ordinals Player 1 has played and $x$ is a element of $\CM$ chosen to elongate the $\bar{a}$-side or the $\bar{b}$-side. Player 2 then must choose $y \in \CM$ to elongate the side opposite which Player 1 has chosen. By the wellfoundedness of the class of ordinals, Player 1 must eventually play the ordinal $0$. After Player 2 responds, the game ends. One says that Player 2 wins this game if and only if the mapping $\bar{a}$ to $\bar{b}$ and the sequence of responses in the game form a partial $\SCRL$-isomorphism. Intuitively, Player 1 winning $\mathrm{EF}_{\alpha}^{\CM,\bar{a},\bar{b}}$ indicates that with $\alpha$-degree of flexibility, Player 1 can compel Player 2 to make a move that violates the $\SCRL$-structure of $\CM$. There is a close relationship between Player 2 having a winning strategy in $\mathrm{EF}_\alpha^{\CM,\bar{a},\bar{b}}$ and the back-and-forth relation $\sim_\alpha$.

By a cardinality consideration, the Scott rank of a tuple in $\CM$ is less than $|\CM|^+$, the cardinal successor of $|\CM|$. Thus the Scott rank of $\CM$ is less than $|\CM|^+$. In particular, if $\CM$ is countable, then $\SR(\CM) < \omega_1$, the first uncountable ordinal. Moreover, Nadel showed there is a close relationship between the definability of $\CM$ and the bound on its Scott rank. Nadel (\cite{Scott-Sentences-and-Admissible-Sets}) showed that $\SR(\CM) \leq \omega_1^{\CM} + 1$ when $\CM$ is considered as a subset of $\omega$. Here $\omega_1^\CM$ is the Church-Kleene ordinal relative to $\CM$ which is the least ordinal $\alpha$ which does not have a wellordering coded as a subset of $\omega$ of ordertype $\alpha$ which is (Turing) computable from $\CM$. It is also the minimal ordinal height of an admissible set containing $\CM$. An admissible set containing $\CM$ is simply a transitive set containing $\CM$ satisfying a weak set theory axiom system called Kripke-Platek $(\KP)$ set theory. The reader can consider an admissible set containing $\CM$ as essentially a miniature universe of set theory containing $\CM$. 

Let $\CM$ be a countable $\SCRL$-structure. The back-and-forth process described above can be used to define a countable infinitary $\SCRL$-formula $\psi_\CM$ so that for any countable $\SCRL$-structure $\CN$, $\CN \models \psi_\CM$ if and only if $\CM$ and $\CN$ are $\SCRL$-isomorphic. The rank of the sentence $\psi_\CM$ is closely related to the Scott rank of $\CM$ and $\psi_\CM$ is roughly the conjunction of all the associated distinguishing formulas for all possible pairs of tuples. For more on the classical and effective Scott analysis for countable structures, see \cite{Model-Theory-an-Introduction}, \cite{Scott-Sentences-and-Admissible-Sets}, and \cite{Computable-Structures-Hyperarithmetic-Hierarchy}.

A particular instance of the above is the study of isometries of metric spaces. The natural first order language $\SCRU$ for metric structures consists of two binary relation symbols for each positive rational $q$ whose intended iterpretations are whether two points have distance less than or more than $q$. (See Definition \ref{language of metric spaces}.) 

A general metric space can have arbitrarily large cardinality and Scott rank. The collection of Polish metric spaces form a very interesting class of metric spaces. Polish metric spaces are complete separable metric spacces.  An uncountable Polish metric space $\CN$ must have cardinality $2^{\aleph_0}$. A priori, one has $\SR(\CN) < (2^{\aleph_0})^+$. There is some hope of doing better. If $\CN$ is a Polish metric space, there is a countable submetric space $\CM \subseteq \CN$ whose completion is $\CN$. By the general theory of Scott analysis mentioned above, $\SR(\CM) < \omega_1$ and in fact $\SR(\CM) < \omega_1^{\CM} + 1$ since $\CM$ is a countable structure. In some sense, $\CM$ has full metric information of it own completion $\CN$. A very natural question asked by Fokina, Friedman, Koerwien, and Nies is whether $\CM$ captures the first order metric structure of it own completion $\CN$ well enough to imply that $\SR(\CN)$ is countable. If so, the author asks whether the Polish metric space $\CN$ with countable dense submetric space $\CM$ satisfies the natural analog of Nadel's effective bounding result for countable structures.

\Begin{question}{main question}
(Fokina, Friedman, Koerwien, Nies) Let $\CN$ be a Polish metric space. Is $\SR(\CN) < \omega_1$? 

(Chan) If $\CN$ is a Polish metric space and $\CM$ is a countable dense submetric space of $\CN$, then is $\SR(\CN) \leq \omega_1^\CM + 1$?
\end{question}

It appears that both questions are still open. (See \cite{Scott-Rank-of-Polish-Metric-Spaces} and \cite{Scott-Rank-of-Polish-Metric-Space-Erratum}.) There are partial answers to these questions. Fokina, Friedman, Koerwien, and Nies showed using some results of Gromov that if $\CN$ is a compact Polish metric space then $\SR(\CN) \leq \omega + 1$. (See also Theorem \ref{compact metric space bound omega + 1}.) Doucha \cite{Scott-Rank-of-Polish-Metric-Spaces} showed that although the cardinality of a Polish metric space $\CN$ is $2^{\aleph_0}$, $\SR(\CN)$ is less than or equal to $\omega_1$, the first uncountable ordinal. Thus Question \ref{main question} is reduced to whether it is possible that there is a Polish metric space $\CN$ with $\SR(\CN) = \omega_1$.

The goal of this paper is to extend a positively answer to Question \ref{main question} for larger classes of Polish metric spaces by producing effective countable bounds on Scott rank. This paper will pursue this in the direction of admissibility theory and the Barwise and Jensen theory of infinitary logic in countable admissible fragments. The advantage of this approach is that one produce not only some desired objects but also an entire miniature universe (of a weak set theory $\KP$) containing these objects. One can then perform a variety of arguments internally and externally of this model of $\KP$ and attempt to reflect internal phenomenon to the real world by absoluteness. This approach gives additional insight on the relation between the Scott rank of the Polish metric space $\CN$ and the definability complexity of any of its countable dense submetric space $\CM$. It also seems to have the benefit of simplifying some technical arguments since the miniature universe of $\KP$ set theory can absorb some combinatorics.

Section \ref{basics} provides the basic definitions. The first order language $\SCRU$ of metric spaces, the back-and-forth relations, the Ehrenfeucht-Fra\"iss\'e game, and the notion of Scott ranks of tuples and structures will be defined. 

Section \ref{bounds for compact polish metric spaces} will give a proof of the following result of Fokina Friedman, Koerwien, and Nies :
\\*
\\*\noindent\textbf{Theorem \ref{compact metric space bound omega + 1}.}
\textit{(Fokina, Friedman, Koerwien, Nies) If $\CM$ is a compact Polish metric space, then $\SR(M) \leq \omega + 1$.}
\\*
\\*\indent This result serves as a warmup for the later theorems in the paper. It contains the approximation idea but is simpler than the subsequent theorems since it involve only playing a single game and there are no admissible sets or ordinals of illfounded models of $\KP$ which are externally illfounded. Nies has mentioned to the author that they had originally proved this result using some theory develop by Gromov. The main combinatorial tool for the proof in this paper is to use the K\"onig lemma to produce a compact approximation system (see Definition \ref{compact approximation system}). The K\"onig lemma is the statement that every finitely branching tree has an infinite path. This is a natural combinatorial principle to apply in this setting since the K\"onig lemma is equivalent to the compactness of a certain closed subsets of $\bairespace$, in it usual topology. (For instance, the weak K\"onig lemma is equivalent to the compactness of the Cantor space, $\cantorspace$.)

A countable metric space $\CM$ along with all its distance relations can be identified with a subset of $\omega$. Let $\CC(\CM)$ denote the metric completion of $\CM$. Note that the elements of the completion of $\CM$ are represented by $\CM$-Cauchy sequence which are essentially reals, i.e. elements of $\bairespace$. The Ehrenfeucht-Fra\"iss\'e game on $\CC(\CM)$ requires Player 2 to give perfect responses in the sense that partial isometries need to be produced. Even if Player 1 plays elements of $\CM$ in the Ehrenfeucht-Fra\"iss\'e game on $\CC(\CM)$, Player 2 may need to respond with an element of $\CC(\CM) \setminus \CM$ to maintain the isometry. However allowing the move to be $\CM$-Cauchy sequences makes the game no longer an integer game. This game can not be absorbed into any countable admissible set. To resolve this, Section \ref{games and ranks} defines a new approximation games $G_\alpha^{f,\bar{a},\bar{b}}$ (see Definition \ref{metric game}), where all the moves are ordinals and elements of $\CM$. Instead of playing perfect responses, Player 2 only needs to produce responses whose errors are no more than that prescribed by some $f : \omega \rightarrow \bbQ^+$ which is a computable function that is strictly decreasing and converging to $0$. Using this game, a new rank $R(\bar{a},\bar{b})$ for pairs of tuples $(\bar{a},\bar{b})$ will be defined. 

It will be shown that $\SR(\bar{a},\bar{b}) \leq R(\bar{a},\bar{b})$. Thus bounding $R(\bar{a},\bar{b})$ will suffice to give a bound on $\SR(\bar{a},\bar{b})$. A metric space is said to be proper if and only if every closed ball is compact. It will be shown that in proper Polish metric spaces, if $(\bar{a},\bar{b})$ is a limit of a sequence of points $\langle (\bar{a}_n,\bar{b}_n) : n \in \omega\rangle$ so that for all $n \in \omega$, $R(\bar{a}_n,\bar{b}_n) > \alpha$, then $R(\bar{a},\bar{b}) > \alpha$. The proof of this result requires playing countably infinite many games simultaneously and thinning out to countably infinite many games at each subsequent stage. In contrast, the main argument of \cite{Scott-Rank-of-Polish-Metric-Spaces} involves $\omega_1$-many simultaneous games and requires a thinning to uncountable nonstationary subsets of $\omega_1$ at subsequent stages. 

Section \ref{admissibility} reviews the basics of admissibility and the theory of infinitary logic in countable admissible fragment including the Jensen's model existence theorem and Barwise compactness. The main technical simplification comes from Fact \ref{ill founded model and automorphisms} which states that if there is an illfounded model $\CA$ containing $\CM$ and two pairs of tuples of $\CM$-Cauchy sequences, $\bar{a}$ and  $\bar{b}$, so that $\CA$ thinks that $R(\bar{a},\bar{b})$ is an $\CA$-ordinal which externally $V$ thinks is $\in^\CA$-illfounded, then one can find (in $V$) an autoisometry of the completion $\CC(\CM)$ taking $\bar{a}$ to $\bar{b}$. This is proved by taking Player 2's winning strategy in $\CA$ for the game associated to the ordinal which is externally illfounded and using it to play forever externally in $V$ to produce an autoisometry. Using this result and an application of Jensen's model existence theorem, one can establish Fact \ref{bounds on metric rank of tuple} which asserts that for any pair of tuples $(\bar{a},\bar{b})$ of $\CM$-Cauchy sequences, $\SR(\bar{a},\bar{b}) \leq R(\bar{a},\bar{b}) < \omega_1^{\CM \oplus \bar{a}\oplus\bar{b}}$. This also gives Doucha's result that $\SR(\CC(\CM)) \leq \omega_1$.

Section \ref{main results} contains the two main theorems of the paper. A metric space $\CN$ is rigid if and only if there are no nontrivial autoisometry of $\CN$. 
\\*
\\*\noindent\textbf{Theorem \ref{rigid structure low scott rank}.} \textit{If $\CM$ is a countable metric space so that $\CC(\CM)$ is a rigid metric space, then $\SR(\CC(\CM)) < \omega_1^\CM$.}
\\*
\\*\indent Usually, in applications of the Jensen's model existence theorem, one can establish the consistency of the relevant theory of an appropriate countable admissible fragment by simply using the real universe as a model. For this theorem, one does not a priori know such an object exists in the real world and so one must establish the consistency of the relevant theory by using Barwise compactness.

If one assumes that the completion of a countable metric space is proper, one can prove the Scott rank of the completion is countable and has the analog of Nadel's effective bound:
\\*
\\*\noindent\textbf{Theorem \ref{proper polish space bound on scott rank}.}
\textit{Let $\CM$ be a metric space on $\omega$. Suppose $\CC(\CM)$ is a proper Polish metric space. Then $\SR(\CM) \leq \omega_1^{\CM} + 1$.} 
\\*
\\*\indent These two theorems extend a positive answer to Question \ref{main question} (even the effective form) for the class of rigid Polish metric spaces and proper Polish metric spaces.

Since an early draft of this paper, Nies and Turetsky (\cite{Logic-Blog-2017}) have produced proofs and expanded some of the results here using recursion theoretic methods. 

The techniques used here to analyze the first order Scott analysis of Polish metric spaces differ in flavor considerably from the classical and effective Scott analysis of countable structures of a countable first order language. The usual technique for finding bounds on Scott ranks for countable structures essentially involves looking at the closure ordinal of an appropriate monotone operator on the countable structure. (See the introduction of \cite{Bounds-on-Continuous-Scott-Rank} for some more details.) \cite{Metric-Scott-Analysis} developed the Scott analysis for continuous logic for metric structures. (See \cite{Metric-Scott-Analysis}, \cite{Model-Theory-for-Metric-Structures} and \cite{Bounds-on-Continuous-Scott-Rank} for the notation and more information.) \cite{Bounds-on-Continuous-Scott-Rank} showed that if $\SCRL$ is a recursive language of continuous logic, $\Omega$ is a weak modulus of continuity with recursive code, $\mathcal{D}$ is a countable $\SCRL$-pre-structure, and $\bar{\mathcal{D}}$ is its completion $\SCRL$-structure, then $\SR_\Omega(\bar{D}) \leq \omega_1^{\mathcal{D}}$. (The definition of Scott rank in \cite{Metric-Scott-Analysis} and \cite{Bounds-on-Continuous-Scott-Rank} is slightly different than the definition used in this paper resulting in a bound that differs by $1$. See the introduction in \cite{Bounds-on-Continuous-Scott-Rank} for a brief explanation.) \cite{Bounds-on-Continuous-Scott-Rank} proves an effective bound on the continuous Scott rank depending on the countable dense substructure which is analogous to the effective bound in the classical Scott analysis for countable structure. Moreover in \cite{Bounds-on-Continuous-Scott-Rank}, the bound is obtained as a closure ordinal of a certain monotone operator on the countable dense substructure which is positive $\Sigma$-definable in an appriopriate admissible set; much like the classical case for countable structures. This may suggest that the metric Scott analysis in continuous logic is the correct and fruitful way to generalize the Scott analysis to Polish metric spaces.

The author would like to acknowledge Alexander Kechris and Andr\'e Nies for comments on earlier drafts of this paper. 

\section{Basics}\label{basics}

\Begin{definition}{language of metric spaces}
The language of a metric space, denoted $\SCRU$, is the following: $\SCRU = \{\dot d_q, \dot d^q : q \in \bbQ^+\}$, where for each $q \in \bbQ^+$, $\dot d_q$ and $\dot d^q$ are binary relation symbols.

If $\CM = (M,d)$ is a metric space on the set $M$ with distance function $d$, then $\CM$ is given the canonical $\SCRU$-structure by defining $(\dot d_q)^\CM(x,y) \Leftrightarrow d(x,y) < q$ and $(\dot d^q)^\CM(x,y) \Leftrightarrow d(x,y) > q$. 
\end{definition}

\Begin{fact}{automorphism and bijective isometry}
Let $\CM$ and $\CN$ be two $\SCRU$-structures which are metric spaces. There is a bijective isometry between $\CM$ and $\CN$ if and only if there is a $\SCRU$-isomorphism between $\CM$ and $\CN$.
\end{fact}

\Begin{definition}{scott rank}
Let $\SCRL$ be any countable first order language. Let $\CM$ be a $\SCRL$-structure. For each ordinal $\alpha$, the relation $\sim_\alpha$ is defined on tuples from $M$ of the same length as follows:

Let $\bar{a} = (a_0, ..., a_{p -1})$ and $\bar{b} = (b_0, ..., b_{p - 1})$, where $p \in \omega$.

$\bar{a} \sim_0 \bar{b}$ if and only if the map sending $a_i$ to $b_i$  for all $i < p$ is a partial $\SCRL$-isomorphism.

$\bar{a} \sim_{\alpha + 1} \bar{b}$ if and only if for all $a \in M$, there exists a $b \in M$ so that $\bar{a}a \sim_\alpha \bar{b}b$ and for all $b \in M$, there exists an $a \in M$ so that $\bar{a}a \sim_\alpha \bar{b}b$. 

If $\beta$ is a limit ordinals, then $\bar{a} \sim_\beta \bar{b}$ if and only if for all $\alpha < \beta$, $\bar{a} \sim_\alpha \bar{b}$. 

Define $\SR(\bar{a},\bar{b}) = \min\{\mu \in \text{ON} : \neg(\bar{a} \sim_\mu \bar{b})\}$ if this set is nonempty. Otherwise $\SR(a,b) = \infty$. 

Define $\SR(\bar{a}) = \sup\{\SR(\bar{a},\bar{b}) : \bar{b} \in {}^{|\bar{a}|} M \wedge \SR(\bar{a},\bar{b}) \neq \infty\}$. 

Finally, $\SR(\CM) = \sup\{\SR(\bar{a}) + 1 : \bar{a} \in {}^{<\omega}M\}$. 
\end{definition}

\Begin{definition}{Ehrenfeuct-Fraisse game}
Let $\SCRL$ be some countable first order language. Let $\CM$ be a $\SCRL$-structure. For some $p \in \omega$, let $\bar{a} = (a_0,...,a_{p - 1})$ and $\bar{b} = (b_0,...,b_{p - 1})$ be tuples from $M$. Let $\alpha$ be an ordinal. The Ehrenfeucht-Fra\"iss\'e game $\EF^{\CM, \bar{a},\bar{b}}_\alpha$ is defined as follows:

If $\alpha = 0$, then Player 2 wins if and only if the map $a_i \mapsto b_i$ for each $i < p$ is a partial $\SCRL$-isomorphism.

If $\alpha > 0$, then Player 1 and 2 play the following:
$$\begin{array}{c | c c c c c c c c }
\bar{a} & (\alpha_0, \Gamma_0 = c_0) & {} & {} (\alpha _1, \Gamma_1 = c_1) & {} & ... & (\alpha_{k - 1}, \Gamma_{k - 1} = d_{k - 1}) & {} \\
\hline
\bar{b} & {}  & \Lambda_0 = d_0 & {} & {} \Lambda_1 = d_1 & ... & {} & \Lambda_{k - 1} = d_{k - 1}
\end{array}$$
$\Gamma$ and $\Lambda$ are formally either the symbol $c$ or $d$. If Player 1 lets $\Gamma$ be $c$, then Player 2 must let the next $\Lambda$ be $d$. If $\Gamma$ is $d$, then $\Lambda$ must be $c$. Each $\Gamma_i$ and $\Lambda_i$ are elements of $M$. $\alpha_0 < \alpha$ and for all $i < k - 1$, $\alpha_{i + 1} < \alpha_i$. The game ends when Player 1 plays $\alpha_{k - 1} = 0$ and $\Gamma_{k - 1}$ and Player 2 responds with $\Lambda_{k - 1}$. When the games ends, a sequence $c_0, ..., c_{k - 1}$ and a sequence $d_0, ..., d_{k - 1}$ have been produced. (The sole purpose of the $\Gamma$ and $\Lambda$ notation is to indicate whether player 1 played $c_i$ (left side associated with $\bar{a}$) or $d_i$ (right side associated with $\bar{b}$) and similarly for Player 2.) Player 2 wins if and only if the map $a_i \mapsto b_i$ for $i < p$ and $c_i \mapsto d_i$ for $i < k$ is a partial $\SCRL$-isomorphism.

The above diagram is a sample play: Here, $\Gamma_0 = c_0$, $\Lambda_0 = d_0$, $\Gamma_1 = c_1$, and $\Lambda_1 = d_1$. This means Player 1 plays first on the left side, Player 2 responds on the right side, Player 1 follows with a play on left side again, and Player 2 responds on the right side, and so forth.
\end{definition}

\Begin{fact}{scott rank and EF games}
Let $\SCRL$ be some countable first order language. Let $\CM$ be an $\SCRL$-structure. Let $\bar{a}$ and $\bar{b}$ be tuples from $\CM$ of the same length. $\SR(\bar{a},\bar{b}) > \alpha$ if and only if Player 2 has a winning strategy in $\EF^{\CM,\bar{a},\bar{b}}_\alpha$. 
\end{fact}

\section{Bounds for Compact Polish Metric Spaces}\label{bounds for compact polish metric spaces}
Fokina, Friedman, Koerwien, and Nies, showed using results of Gromov about metric spaces that $\SR(M) \leq \omega + 1$, when $M$ is a compact Polish metric space.

This section will give a proof of this result using K\"onig's lemma. This result will require looking at partial maps that are not isometries but have a predetermined error in distances. This argument is a simple approximation idea using a single game which will be a warmup for the later results on proper metric spaces that combines the approximation idea with admissibility, ordinals of admissible sets which are externally illfounded, and infinitely many games.

\Begin{definition}{compact approximation system}
Let $M$ be a compact Polish metric space. Let $p \in \omega$ and $\bar{a} = (a_0, ..., a_{p - 1})$ and $\bar{b} = (b_0, ..., b_{p - 1})$ be tuples of elements from $M$. Let $(A_n : n \in \omega)$ be a sequence of finite subsets of $M$ with the property that for all $n \in \omega$, $A_n \subseteq A_{n + 1}$ and $\bigcup_{z \in A_n} B_{2^{-n}}(z) = M$. 

A compact approximation system (for $M$, $\bar{a}$, and $\bar{b}$ with respect to $(A_n : n \in \omega)$) is a sequence $(\varphi_n : n \in \omega)$ with the following properties:

(i) $\varphi_n : A_n \rightarrow A_n$. 

(ii) For all $i < p$ and $z \in A_n$, $|d(a_i, z) - d(b_i, \varphi_n(z))| < 2^{-n}$. 

(iii) For all $m \leq n$, for all $y \in A_m$ and $z \in A_n$, 
$$|d(y,z) - d(\varphi_m(y), \varphi_n(z)| < 2^{-m} + 2^{-n}.$$

(iv) For all $n \in \omega$, $\bigcup_{z \in A_n} B_{2^{-(n - 1)}}(\varphi_n(z)) = M$. 

\noindent A $k$-compact approximation system is a sequence $(\varphi_n : n \leq k)$ satisfying the above properties below $k$.
\end{definition}

\Begin{lemma}{compact approximation system yield autoisometry}
Suppose $(\varphi_n : n \in \omega)$ is a compact approximation system for $M$, $\bar{a}$, and $\bar{b}$ with respect to $(A_n : n \in \omega)$, then there is an autoisometry $\Phi : M \rightarrow M$ such that for all $i < p$, $\Phi(a_i) = b_i$. 
\end{lemma}

\begin{proof}
Let $x \in M$. Let $(x_n : n \in \omega)$ be a sequence with the property that for all $n \in \omega$, $x_n \in A_n$ and $\lim_{n \rightarrow \infty} x_n = x$. Define $\Phi(x) = \lim_{n \rightarrow \infty} \varphi_n(x_n)$. 

It remains to show that $\Phi$ is well-defined and $\Phi$ is an autoisometry with $\Phi(\bar{a}) = \bar{b}$. 

First to show that $(\varphi_n(x_n) : n \in \omega)$ is a Cauchy sequence: Let $m \leq n$. By (iii)
$$d(\varphi_m(x_m), \varphi_n(x_n)) < d(x_m,x_n) + 2^{-m} + 2^{-n}$$
$$\leq d(x_m, x) + d(x_n, x) + 2^{-m} + 2^{-n}$$
Since $\lim_{n \rightarrow \infty} x_n = x$, this shows that $(\varphi_n(x_n) : n \in \omega)$ is a Cauchy sequence.

Next, to show that $\Phi(x)$ is independent of the sequence $(x_n : n \in \omega)$ which is used to define it: Suppose $(y_n : n \in \omega)$ is another sequence with the property that $y_n \in A_n$ and $\lim_{n \rightarrow \infty} y_n = x$. Since $x_n,y_n \in A_n$, (iii) states
$$|d(\varphi_n(x_n), \varphi_n(y_n)) - d(x_n,y_n)| < 2^{-n} + 2^{-n}$$
Therefore
$$d(\varphi_n(x_n), \varphi_n(y_n)) < 2^{-(n - 1)} + d(x_n,y_n)$$
$$\leq 2^{-(n - 1)} + d(x_n,x) + d(x, y_n)$$
Since $\lim_{n \rightarrow \infty} x_n = \lim_{n \rightarrow \infty} y_n = x$, the above shows that $\lim_{n \rightarrow \infty} d(\varphi_n(x_n),\varphi_n(y_n)) = 0$. So $\lim_{n \rightarrow \infty} \varphi_n(x_n) = \lim_{n \rightarrow \infty} \varphi_n(y_n)$. This shows that $\Phi$ is a well-defined function.

Next, to show that for any $i < p$, $\Phi(a_i) = b_i$: Let $(x_n : n \in \omega)$ be a sequence such that $x_n \in A_n$ and $\lim_{n \rightarrow \infty} x_n = a_i$. By (ii),
$$d(\varphi_n(x_n), b_i) < d(x_n,a_i) + 2^{-n}$$ 
Since $\lim_{n \rightarrow \infty} x_n = a_i$, this shows that $\Phi(a_i) = \lim_{n \rightarrow \infty} \varphi_n(x_n) = b_i$. 

Next, to show that $\Phi$ is an isometry: Suppose $\lim_{n \rightarrow \infty} e_n = e$ and $\lim_{n \rightarrow \infty} f_n = f$. Then
$$|d(e,f) - d(e_n,f_n)| = |d(e,f) - d(e,f_n) + d(e,f_n) - d(e_n,f_n)|$$
$$\leq |d(e,f) - d(e,f_n)| + |d(e,f_n) - d(e_n, f_n)| \leq d(f,f_n) + d(e,e_n)$$
Therefore, $\lim_{n \rightarrow \infty} d(e_n,f_n) = d(e,f)$. 

Now suppose $x,y \in M$. Let $(x_n : n \in \omega)$ and $(y_n : n \in \omega)$ be such that $x_n,y_n \in A_n$ and $\lim_{n \rightarrow \infty} x_n = x$ and $\lim_{n \rightarrow \infty} y_n = y$. By (iii),
$$|d(\varphi_n(x_n), \varphi_n(y_n)) - d(x_n,y_n)| < 2^{-(n - 1)}$$
which implies that $\lim_{n \rightarrow \infty}d(\varphi_n(x_n),\varphi_n(y_n)) = \lim_{n \rightarrow \infty} d(x_n,y_n)$. So using the observation of the previous paragraph for the first and third equality,
$$d(\Phi(x),\Phi(y)) = \lim_{n \rightarrow \infty} d(\varphi_n(x_n), \varphi_n(y_n)) = \lim_{n \rightarrow \infty} d(x_n,y_n) = d(x,y).$$ 
This shows that $\Phi$ is an isometry.

Finally to show that $\Phi$ is surjective: Let $y \in M$. By (iv), for all $n \in \omega$, $\bigcup_{z \in A_n} B_{2^{-(n - 1)}}(\varphi_n(z)) = M$. For each $n \in \omega$, let $x_n \in A_n$ be such that $d(y, \varphi_n(x_n)) < 2^{-(n - 1)}$. Observe that $(x_n : n \in \omega)$ is a Cauchy sequence. To see this, by (iii),
$$d(x_m,x_n) < 2^{-m} + 2^{-n} + d(\varphi_m(x_m), \varphi_n(x_n))$$
$$\leq 2^{-m} + 2^{-n} + d(\varphi_m(x_m), y) + d(y, \varphi_n(x_n))$$
$$\leq 3(2^{-m}) + 3(2^{-n})$$
Therefore, let $x = \lim_{n \rightarrow \infty} x_n$. Then $y = \lim_{n \rightarrow \infty} \varphi_n(x_n) = \Phi(x)$. This shows that $\Phi$ is surjective and completes the proof of the lemma.
\end{proof}

\Begin{lemma}{SR condition implies compact approximation system}
Let $M$ be a compact Polish metric space. Let $\bar{a} = (a_0, ..., a_{p - 1})$ and $\bar{b} = (b_0, ..., b_{p - 1})$ be tuples from $M$. Let $(A_n : n \in \omega)$ be a sequence of finite subsets of $M$ so that for all $n \in \omega$, $A_n \subseteq A_{n + 1}$ and $\bigcup_{z \in A_n} B_{2^{-n}}(z) = M$. Suppose $\SR(\bar{a},\bar{b}) > \omega$, then there is a compact approximation system for $M$, $\bar{a}$, $\bar{b}$ with respect to $(A_n : n \in \omega)$. 
\end{lemma}

\begin{proof}
Define $J$ to be the tree of all $k$-compact approximation system for $M$, $\bar{a}$, and $\bar{b}$ with respect to $(A_n : n \in \omega)$, where $k$ varies over $\omega$. $J$ is ordered by $(\sigma_i : i \leq m) \preceq_J (\tau_i : i \leq n)$ if and only if $m \leq n$ and for all $i \leq m$, $\sigma_i = \tau_i$. As each $A_n$ is finite, $J$ is a finitely branching tree. Any infinite path through $J$ would be a compact approximation system. By K\"onig's lemma, $J$ would have an infinite path if $J$ was infinite. 

As $\SR(\bar{a},\bar{b}) > \omega$, fix a winning strategy for Player 2 in $\EF^{M,\bar{a},\bar{b}}_\omega$. To show $J$ is infinite, it suffices to show that there is a $k$-compact approximation system for each $k \in \omega$. Let $L = |A_k|$. Enumerate $A_k = \{c_i : i < L\}$. Consider the following game of $\EF^{M,\bar{a},\bar{b}}_\omega$ where Player 1 plays $(L - i, c_i)$ (i.e. Player 1 plays $\Gamma = c$) and Player 2 always responds with the winning strategy:
$$\begin{array}{c | c  c c c c c c c}
\bar{a} & (L, c_{0}) & {} & (L - 1, c_{1}) & {} & ... & (1, c_{L - 1}) & {} \\
\hline
\bar{b} & {} & d_{0} & {} & d_{1} & ... & {} & d_{L - 1}
\end{array}$$
(Note that the last ordinal played is 1, which allows player 1 to play one more time.)

(Recall that if $n \leq k$, $A_n \subseteq A_k$.) For each $n \leq k$ and $c_i \in A_n$, define $\varphi_n(c_i)$ to be some element of $A_n$ so that $d(d_i, \varphi_n(c_i)) < 2^{-n}$, which is possible since $\bigcup_{z \in A_n}B_{2^{-n}}(z) = M$. This completes the definition of $(\varphi_n : n \leq k)$.

Now to check that $(\varphi_n : n \leq k)$ is a $k$-compact approximation system: (i) is clearly true. For (ii): Pick some $i < p$ and $c_j \in A_n$,
$$|d(a_i, c_j) - d(b_i, \varphi_n(c_j))| = |d(a_i, c_j) - d(b_i,d_j) + d(b_i, d_j) - d(b_i, \varphi_n(c_j))|$$
$$\leq |d(a_i,c_j) - d(b_i,d_j)| + |d(b_i,d_j) - d(b_i,\varphi_n(c_j))|$$
$$\leq 0 + d(d_j,\varphi_n(c_j)) < 2^{-n}$$
since Player 2 used its winning strategy for $\mathrm{EF}^{M,\bar{a},\bar{b}}_\omega$ and by the definition of $\varphi_n(c_j)$.

For (iii): Let $m \leq n \leq k$, $c_i \in A_m$, and $c_j \in A_n$. 
$$|d(c_i,c_j) - d(\varphi_m(c_i), \varphi_n(c_j))| = |d(c_i,c_j) - d(d_i,d_j) + d(d_i,d_j) - d(\varphi_m(c_i), \varphi_n(c_j))|$$
$$\leq |d(c_i,c_j) - d(d_i,d_j)| + |d(d_i,d_j) + d(\varphi_m(c_i), \varphi_n(c_j))|$$
$$= 0 + |d(d_i,d_j) - d(\varphi_m(c_i),\varphi_n(c_j))|$$
$$= |d(d_i,d_j) - d(\varphi_m(c_i), d_j) + d(\varphi_m(c_i),d_j) - d(\varphi_m(c_i),\varphi_n(c_j))|$$
$$\leq |d(d_i,d_j) - d(\varphi_m(c_i), d_j)| + |d(\varphi_m(c_i), d_j) - d(\varphi_m(c_i), \varphi_n(c_j))|$$
$$\leq d(d_i,\varphi_m(c_i)) + d(d_j, \varphi_n(c_j)) < 2^{-m} + 2^{-n}$$

For (iv): Fix $n \leq k$. Let $R = |A_n|$. Let $(c_{i_l} : l < R)$ be the subsequence enumerating $A_n \subseteq A_k$. Suppose that $\bigcup_{l < R} B_{2^{-(n - 1)}}(\varphi_n(c_{i_l})) \subsetneq M$. Then there is some $y$ so that for all $l < R$, $d(y, \varphi_n(c_{i_l})) \geq 2^{-(n - 1)}$. 

Note that for all $l < R$, $d(d_{i_l}, y) \geq 2^{-n}$: To see this, suppose that there were some $l < R$ so that $d(d_{i_l}, y) < 2^{-n}$. Then
$$d(\varphi_n(c_{i_l}), y) \leq d(\varphi_n(c_{i_l}), d_{i_l}) + d(d_{i_l}, y) < 2^{-n} + 2^{-n} = 2^{-(n - 1)}$$
This is a contradiction.

Let $d_L = y$. Now continue playing the game $\EF^{M,\bar{a},\bar{b}}_\omega$ one more time as follows:
$$\begin{array}{c | c  c c c c c c c c}
\bar{a} & (L, c_{0}) & {} & (L - 1, c_{1}) & {} & ... & (1, c_{L - 1}) & {} & {} & c_L \\
\hline
\bar{b} & {} & d_{0} & {} & d_{1} & ... & {} & d_{L - 1} & (0,d_L) & {} 
\end{array}$$
(This means that $\Gamma = d$ in the last time that Player 1 moves, i.e. Player 1 played on the right side.) Let $c_L$ be the response by Player 2 using its winning strategy. The claim is that the map induced by this play is not a partial isometry. To see this: Since $\bigcup_{l < R} B_{2^{-n}}(c_{i_l}) = M$, there is some $l < R$ so that $d(c_L, c_{i_l}) < 2^{-n}$. Then $d(y, d_{i_l}) = d(d_L, d_{i_l}) < 2^{-n}$. This contradicts the result of the previous paragraph. 

This completes the proof of the lemma.
\end{proof}

As an immediate corollary, one obtains the result of Fokina, Friedman, Koerwien, and Nies on Scott ranks of compact Polish metric spaces.

\Begin{theorem}{compact metric space bound omega + 1}
(Fokina, Friedman, Koerwien, and Nies) If $M$ is a compact Polish metric space, then $\SR(M) \leq \omega + 1$. 
\end{theorem}

The next few sections will be concerned with finding an effective bound on the Scott rank of proper Polish metric spaces.

\section{Games and Ranks}\label{games and ranks}
For the rest of the paper, let $\CM$ be a countably infinite metric space. By taking a bijection, one may assume that the domain of the metric space $\CM$ is $\omega$. By considering the domain of $\CM$ as $\omega$, $\CM$ can be coded as a real, i.e. an element of $\bairespace$, by coding all the interpretations of symbols of $\SCRU$ as relations on $\omega$ in some fixed way. This section will consider $\CM$ as a metric space; however, in the following section, one will occasionally refer to $\CM$ as a real which codes the structure in the above way.

If $\CM$ is a metric space, then $\CC(\CM)$ denotes the metric space completion of $\CM$.

Based on the method of constructing a bijective isometry in \cite{Scott-Rank-of-Polish-Metric-Spaces} Lemma 2.3, one defines a new rank  on tuples of elements of $\CC(\CM)$ depending on whether Player 2 has a winning strategy in some game on $\CM$ (essentially on $\omega$). Since the Ehrenfeucht-Fra\"{i}sse game requires the construction of partial $\SCRU$-isomorphisms, even if Player 1 always plays elements of $\CM$, Player 2 generally needs to respond with a $\CM$-Cauchy sequence (essentially an element of $\bairespace$). This makes the definability and absoluteness property of Scott rank (in respect to descriptive set theoretic complexity) quite difficult to determine. A priori, it seems quite possible that playing the Ehrenfeucht-Fra\"{i}se game in different models of set theory with either more or less Cauchy sequences could affect the outcome of the game. This new game will be played on $\CM$ so any model of set theory containing $\CM$ (which includes the interpretations of the symbols of $\SCRU$) will play these games correctly. The Ehrenfeucht-Fra\"{i}sse game asks Player 1 to play perfectly in the sense that it must produce partial isomorphisms; this new game will be ostensibly easier for Player 2 since it demands only the response be appropriately close to Player 1's move.

The following convention in variable naming will be used: The variables $a$, $c$, and $x$ will denote objects on the left side. The variable $b$, $d$ and $y$ will denote objects on the right side. Throughout the paper there may be bars or subscripts attached to these variables but they will always denote plays on the sides indicated.

\Begin{definition}{metric game}
Let $\CM$ be a metric space on $\omega$. Let $\CC(\CM)$ be its completion. Let $\bar{a} = (a_0, ..., a_{p - 1})$ and $\bar{b} = (b_0, b_1, ..., b_{p - 1})$ be tuples of elements of $\CC(\CM)$. Let $f : \omega \rightarrow \bbQ^+$ denote a recursive (i.e. computable) strictly decreasing function converging to $0$. Let $\alpha$ be an ordinal.

Define the following game $G^{f,\bar{a},\bar{b}}_\alpha$:
$$\begin{array}{c | c c c c c c c c c c c c c }
\bar{a} & (\alpha_0,c_0) & {} & {} & c_1 & (\alpha_2,c_2) & {} & {} & c_3 & \dots & (\alpha_{k - 1},c_{k - 1}) & {}  \\
\hline
\bar{b} & {} & d_0 & (\alpha_1, d_1) & {} & {} & d_2 & (\alpha_3, d_3) & {} & \dots  & {} & d_{k - 1}
\end{array}
$$
Player 1 and Player 2 alternatingly play $(\alpha_0,c_0)$, $d_0$, $(\alpha_1, d_1)$, $c_1$, $(\alpha_2, c_2)$, $d_2$, $(\alpha_3, d_3)$, $c_3$, ..., $(\alpha_{k - 1}, r_{k - 1})$, $s_{k - 1}$, where $r = c$ and $s = d$ if $k$ is odd and $r = d$ and $s = c$ if $k$ is even. For all $i < k$, $\alpha_i$ is an ordinal less than $\alpha$. For all $i < k$, $c_i$ and $d_i$ are elements of $\CM$. Since $\CM$ is a metric space on $\omega$, $c_i$ and $d_i$ are natural numbers. For all $i < k - 1$, $\alpha_{i + 1} < \alpha_{i}$. The game ends when Player 1 plays $\alpha_{k -1} = 0$ and Player 2 responds. Player 2 wins if and only if the following holds:

(I) For all $i < p$ and $j < k$, $|d(a_i, c_j) - d(b_i, d_j)| < f(j)$.

(II) For all $i,j < k$, $|d(c_i,c_j) - d(d_i,d_j)| < f(i) + f(j)$. 

\noindent In the above, the distance function $d$ refers to the distance function of $\CC(\CM)$. 
\end{definition}

\Begin{definition}{metric rank}
Let $\CM$ be a metric space on $\omega$. Let $\mathrm{REC}$ be the set of all $f : \omega \rightarrow \bbQ^+$ which are recursive, strictly decreasing, and converge to $0$. 

Let $\bar{a}$ and $\bar{b}$ be two tuples of elements from $\CC(\CM)$ of the same length. Let $f \in \mathrm{REC}$. Say $\bar{a} \sim^f_\alpha \bar{b}$ if and only if Player 2 has a winning strategy in $G_\alpha^{f, \bar{a},\bar{b}}$. 

Define 
$$\RR(\bar{a},\bar{b}) = \min \{\mu \in \text{ON} : (\exists f \in \text{REC})\neg(\bar{a} \sim^f_\mu \bar{b})\}$$
if the above set is nonempty. Otherwise let $\RR(\bar{a},\bar{b}) = \infty$.
\end{definition}

The use of the class $\mathrm{REC}$ of recursive, strictly decreasing functions taking values in the positive rational numbers and converging to $0$ is merely for convenience. The important property is that these functions used in the following sections is that they are coded in any admissible set. By inspecting the proof, one can find a much smaller class of such functions that would be adequate for the following arguments.

Next, the relationship between $\RR$ and $\SR$ will be determined:

\Begin{fact}{relation between scott rank and metric rank}
Let $\CM$ be a metric space on $\omega$. Let $\bar{a}$ and $\bar{b}$ be two tuples of elements of $\CC(\CM)$ of the same length. For all ordinals $\alpha$ and $f \in \mathrm{REC}$, if $\bar{a} \sim_\alpha \bar{b}$, then $\bar{a} \sim_\alpha^f \bar{b}$. Hence $\SR(\bar{a},\bar{b}) \leq \RR(\bar{a},\bar{b})$. 
\end{fact}

\begin{proof}
Let $\alpha < \SR(\bar{a},\bar{b})$. Let $f \in \mathrm{REC}$. A winning strategy for Player 2 in $G^{f,\bar{a},\bar{b}}_\alpha$ will be produced. 

The idea is that $\bar{a} \sim_\alpha \bar{b}$ allows Player 2 to find perfect responses (albeit in $\CC(\CM)$ not $\CM$) in the Ehrenfeucht-Fra\"{i}sse game in the sense that the responses form a partial isometry. So Player 2 will respond in the game $G^{f,\bar{a},\bar{b}}_\alpha$ by simply choosing some element of $M$ (i.e. $\omega$) which is sufficiently close to the perfect response given by $\bar{a} \sim_\alpha \bar{b}$. The details follows:

Consider a play of $G_\alpha^{f,\bar{a},\bar{b}}$. Now suppose $(\alpha_0,c_0)$, $y_0$, $d_0$, $(\alpha_1, d_1)$, $x_1$, $c_1$, ..., $(\alpha_{j - 1}, d_{j - 1})$, $x_{j - 1}$, and $c_{j - 1}$ has appeared in the construction thus far (assuming $j$ is even) and satisfies the following: For all $i < j$, 

(i) $\bar{a}\hat{\ }c_0\hat{\ }x_1 \hat{\ } c_2 \hat{\ } x_3 \hat{\ } ...\hat{\ } c_{i} \sim_{\alpha_i} \bar{b} \hat{\ } y_0\hat{\ } d_1 \hat{\ } y_2 \hat{\ } d_3 \hat{\ } ... \hat{\ } y_i$ if $i$ is even. If if $i$ is odd, the same holds with the last $c_i$ replaced by $x_i$ and $y_i$ replaced by $d_i$. 

(ii) $d(y_i, d_i) < f(i)$ if $i$ is even or $d(x_i,c_i) < f(i)$ if $i$ is odd.

Assuming that $\alpha_{j - 1} \neq 0$, suppose Player 1 chooses to play $(\alpha_j, c_j)$ where $\alpha_j < \alpha_{j - 1}$. Since 
$$\bar{a}\hat{\ }c_0\hat{\ }x_1 \hat{\ } c_2 \hat{\ } x_3 \hat{\ } ...\hat{\ } y_{j - 1} \sim_{\alpha_{j - 1}} \bar{b} \hat{\ } y_0\hat{\ } d_1 \hat{\ } y_2 \hat{\ } d_3 \hat{\ } ... \hat{\ } d_{j - 1}$$
one can find some $y_j$ so that 
$$\bar{a}\hat{\ }c_0\hat{\ }x_1 \hat{\ } c_2 \hat{\ } x_3 \hat{\ } ...\hat{\ } y_{j - 1} \hat{\ } c_j \sim_{\alpha_j} \bar{b} \hat{\ } y_0\hat{\ } d_1 \hat{\ } y_2 \hat{\ } d_3 \hat{\ } ... \hat{\ } y_{j - 1} \hat{\ } y_j$$
Now let $d_j \in M$ be chosen so that $d(y_j, d_j) < f(j)$. 

Now continue this process as long as Player 1 has not played the ordinal $0$. Of course, depending on whether the stage is even or odd, the variable needs to be appropriately changed.

At some stage $k$, Player 1 will have played $\alpha_{k - 1} = 0$. After Player 2 responds, the process finishes.

The claim is that the following play of $G^{f,\bar{a},\bar{b}}_\alpha$ is winning for Player 2:
$$\begin{array}{c | c c c c c c c c c c c c c }
\bar{a} & (\alpha_0,c_0) & {} & {} & c_1 & (\alpha_2,c_2) & {} & {} & c_3 & \dots & (\alpha_{k - 1},c_{k - 1}) & {}  \\
\hline
\bar{b} & {} & d_0 & (\alpha_1, d_1) & {} & {} & d_2 & (\alpha_3, d_3) & {} & \dots  & {} & d_{k - 1}
\end{array}
$$

First pick some $i < p$ and $j < k$. Without loss of generality, suppose $j$ is even. Then
$$|d(a_i,c_j) - d(b_i,d_j)| = |d(a_i,c_j) - d(b_i,y_j) + d(b_i,y_j) - d(b_i,d_j)|$$
Recall $d(a_i,c_j) = d(b_i,y_j)$ hence
$$= |d(b_i,y_j) - d(b_i, d_j)| \leq d(y_j, d_j) < f(j)$$

Now pick some $i,j < k$. Without loss of generality suppose $i$ is even and  $j$ is odd. Then
$$|d(c_i,c_j) - d(d_i,d_j)| = |d(c_i,c_j) - d(c_i,x_j) + d(c_i,x_j) - d(d_i,d_j)|$$
Recall that $d(c_i,x_j) = d(y_i, d_j)$. Therefore
$$= |d(c_i,c_j) - d(c_i,x_j) + d(y_i,d_j) - d(d_i,d_j)|$$
$$\leq |d(c_i,c_j) - d(c_i,x_j)| + |d(y_i,d_j) - d(d_i,d_j)|$$
$$\leq d(c_j,x_j) + d(y_i, d_i) < f(j) + f(i)$$
All the other even and odd combinations are handled similarly.
\end{proof}

\Begin{definition}{proper metric spaces}
A metric space $\CN$ is \textit{proper} if and only if if $\{y : d(x,y) \leq r\}$ is compact for all $x \in N$ and $r \in \bbR$.
\end{definition}

The key property of proper metric spaces that will be used is that every bounded sequence has a convergent subsequence.

In the following, suppose $\CN$ is some metric space with distance $d_\CN$. Then the distance on ${}^k\CN$ is defined as $d_{{}^k\CN}(\bar{a},\bar{b}) = \sum_{i = 0}^{k - 1} d_\CN(a_i,b_i)$. 

Now suppose $\bar{a}$ and $\bar{b}$ are two tuples of length $p$ of elements of $\CC(\CM)$.  The next technical lemma asserts that if $\alpha$ is an ordinal, $\CC(\CM)$ is a proper metric space, and $(\bar{a},\bar{b})$ is a limit (in the $(\CC(\CM))^{2p}$ metric) of points of the form $(\bar{e},\bar{f})$ so that $\RR(\bar{e},\bar{f}) > \alpha$, then  $\RR(\bar{a},\bar{b}) > \alpha$. 

\Begin{fact}{proper polish space limit certain rank implies certain rank}
Let $\CM$ be a metric space on $\omega$. Suppose $\CC(\CM)$ is a proper metric space. Let $\alpha$ be an ordinal. Let $\bar{a}$ and $\bar{b}$ be two tuples of elements of $\CC(\CM)$ of the same length $p$. Suppose that $(\bar{a},\bar{b})$ is the limit of the sequence $\langle (\bar{a}_n,\bar{b}_n) : n \in \omega \rangle$ in $(\CC(\CM))^{2p}$ so that for all $n \in \omega$, $\RR(\bar{a}_n,\bar{b}_n) > \alpha$. Then $\RR(\bar{a},\bar{b}) > \alpha$. 
\end{fact}

\begin{proof}
Fix $f \in \mathrm{REC}$. Let $g \in \mathrm{REC}$ be defined by $g(n) = \frac{f(n)}{3}$ for all $n \in \omega$. 

Suppose $\bar{a}$ and $\bar{b}$ take the following form: $\bar{a} = (a_0, ..., a_{p -1})$ and $\bar{b} = (b_0, ..., b_{p - 1})$. For each $n \in \omega$, suppose $\bar{a}_n$ and $\bar{b}_n$ take the form: $\bar{a}_n = (a_0^n, ..., a_{p - 1}^n)$ and $\bar{b}_n = (b_0^n,...,b_{p - 1}^n)$. 

It suffices to show that for all ordinals $\alpha$, if $\bar{a}_n \sim_\alpha^g \bar{b}_n$ for all $n \in \omega$, then $\bar{a} \sim_\alpha^f \bar{b}$. 

Fix a winning strategy for Player 2 in each game $G^{g,\bar{a}_n,\bar{b}_n}_\alpha$. In the following proof, when a response from Player 2 in $G^{g,\bar{a}_n,\bar{b}_n}_\alpha$ is required, it is always assumed it is taken from this fixed winning strategy. 

Now a winning strategy for $G^{f,\bar{a},\bar{b}}_\alpha$ will be described:

By refining $\langle (\bar{a}_n,\bar{b}_n) : n \in \omega\rangle$ to a subsequence, one may assume that 
$$d_p((\bar{a}_n,\bar{b}_n), (\bar{a},\bar{b})) < \frac{1}{n} \ \ \ (\star)$$
where $d_p$ is the metric on $\CC(\CM)^{2p}$ mentioned above which is defined by summing the distance in each coordinate.

Let $A_{-1} = \omega$. 

Fix $j$ and suppose the following have been constructed: For all $i < j$, $A_i$ and $\alpha_i$ have been defined. If $i$ is even, then $c_i$, $y_n$, and $d_i^n$ for each $n \in A_i$ have been constructed. If $i$ is odd, then $d_i$, $x_n$, and $c_i^n$ for each $n \in A_i$ have been constructed. These objects satisfy the following:

(i) For all $i < j$, $\min A_i > \frac{1}{g(i)}$ and $A_i$ is an infinite subset of $\omega$. Hence by property $(\star)$ on the sequence, one has $d((\bar{a}_n,\bar{b}_n),(\bar{a},\bar{b})) < g(i)$, for all $n \in A_i$. 

(ii) For all $i < j - 1$, $A_{i + 1} \subseteq A_i$. 

(iii) If $i$ is even, for all $n \in A_i$, $d(y_i, d_i^n) < g(i)$. If $i$ is odd, for all $n \in A_i$, $d(x_i, c_i^n) < g(i)$. 

(iv) For each $i < j$ and $n \in A_i$, the following is a play in $G^{g,\bar{a}_n,\bar{b}_n}_\alpha$ according to the fixed winning strategy for Player 2:
$$\begin{array}{c | c c c c c c c c c c c c c }
\bar{a}_n & (\alpha_0,c_0) & {} & {} & c_1^n & (\alpha_2,c_2) & {} & {} & c_3^n & \dots & (\alpha_{i},c_{i}) & {}  \\
\hline
\bar{b}_n & {} & d_0^n & (\alpha_1, d_1) & {} & {} & d_2^n & (\alpha_3, d_3) & {} & \dots  & {} & d_{i}^n
\end{array}
$$
when $i$ is even. A similar diagram when $i$ is odd with the appropriate variable change.

Without loss of generality, suppose that $j$ is odd. Suppose that $\alpha_{j - 1} \neq 0$. Suppose Player 1 plays $(\alpha_j, d_j)$ where $\alpha_j < \alpha_{j - 1}$. 

For each $n \in A_{j - 1} \cap (\frac{1}{g(j)}, \infty)$, let $c_j^n$ be the response of Player 2 in the following play of $G^{g,\bar{a}_n,\bar{b}_n}_\alpha$ according to the fixed winning strategy:
$$\begin{array}{c | c c c c c c c c c c c c c }
\bar{a}_n & (\alpha_0,c_0) & {} & {} & c_1^n & (\alpha_2,c_2) & {} & {} & c_3^n & \dots & (\alpha_{j - 1},c_{j - 1}) & {} & {} & c_j^n  \\
\hline
\bar{b}_n & {} & d_0^n & (\alpha_1, d_1) & {} & {} & d_2^n & (\alpha_3, d_3) & {} & \dots  & {} & d_{j - 1}^n & (\alpha_j, d_j) & {}
\end{array}
$$
Claim: $L_j = \{c_j^n : n \in A_{j - 1} \cap (\frac{1}{g(j)},\infty)\}$ is a bounded set 

To see this: Note that for all $n \in A_{j - 1}$, 
$$d(a_0, c_j^n) \leq d(a_0, a_0^n) + d(a_0^n, c_j^n)$$
Since Player 2 plays according to its winning strategy in $G^{g,\bar{a}_n,\bar{b}_n}$, $|d(a_0^n, c_j^n) - d(b_0^n, d_j)| < g(j)$ so
$$ \leq d(a_0, a_0^n) + d(b_0^n, d_j) + g(j) \leq d(a_0,a_0^n) + d(b_0^n, b_0) + d(b_0, d_j) + g(j)$$
Since $n > \frac{1}{g(j)}$ and property $(\star)$ on the sequence, $d(a_0, a_0^n) < g(j)$, $d(b_0, b_0^n) < g(j)$. Thus
$$< 3 g(j) + d(b_0,d_j) = 3 \frac{f(j)}{3} + d(b_0, d_j) = f(j) + d(b_0,d_j)$$
This shows that for all $n \in A_{j - 1}$, $d(a_0,c_j^n) < f(j) + d(b_0,d_j)$. Hence $L_j$ is bounded. The claim has been established.

Since $L_j$ is bounded and $\CC(\CM)$ is a proper metric space, the sequence $\langle c_j^n : n \in A_{j - 1} \cap (\frac{1}{g(j)}, \infty)\rangle$ has a convergent subsequence. Let $x_j \in \CC(\CM)$ be a limit point of a convergent subsequence. Let $A_{j} \subseteq A_{j - 1} \cap (\frac{1}{g(j)}, \infty)$ be an infinite set so that $d(x_j, c_j^n) < g(j)$ for all $n \in A_j$. This completes the recursive construction at stage $j$.

Continue this construction until at some point Player 1 plays $\alpha_{k - 1} = 0$. 

Now the claim is that the following is a winning play for Player 2 in $G^{f,\bar{a},\bar{b}}_\alpha$:
$$\begin{array}{c | c c c c c c c c c c c c c }
\bar{a} & (\alpha_0,c_0) & {} & {} & c_1^{\min A_1} & (\alpha_2,c_2) & {} & {} & c_3^{\min A_3} & \dots & (\alpha_{k - 1},c_{k - 1}) & {}  \\
\hline
\bar{b} & {} & d_0^{\min A_0} & (\alpha_1, d_1) & {} & {} & d_2^{\min A_2} & (\alpha_3, d_3) & {} & \dots  & {} & d_{k - 1}^{\min A_{k - 1}}
\end{array}
$$
Let $l < p$ and $i < k$. Without loss of generality, suppose $i$ is even:
$$|d(a_l, c_i) - d(b_l, d_i^{\min A_i})|$$
$$= |d(a_l, c_i) - d(a_l^{\min A_i}, c_i) + d(a^{\min A_i}_l, c_i) - d(b^{\min A_i}_l, d_i^{\min A_i}) + d(b^{\min A_i}_l, d_i^{\min A_i}) - d(b_l,d_i^{\min A_i})|$$
$$\leq |d(a_l, c_i) - d(a^{\min A_i}_l, c_i)| + |d(a^{\min A_i}_l, c_i) - d(b^{\min A_i}_l, d_i^{\min A_i})| + |d(b^{\min A_i}_l, d_i^{\min A_i}) - d(b_l,d_i^{\min A_i})|$$
$$\leq d(a_l, a^{\min A_i}_l) + |d(a^{\min A_i}_l, c_i) - d(b^{\min A_i}_l, d_i^{\min A_i})| + d(b_{l}^{\min A_i}, b_l)$$
The first and third terms are less than $g(i)$ by (i). The middle term is less that $g(i)$ since these are responses that come from the winning strategy of $G^{g,\bar{a}_{\min A_i},\bar{b}_{\min A_i}}_\alpha$.
$$\leq g(i) + g(i) + g(i) = 3 g(i) = 3 \frac{f(i)}{3} = f(i)$$

Now let $i < j < k$. As an example, assume $i$ is even and $j$ is odd:
$$|d(c_i, c_j^{\min A_j}) - d(d_i^{\min A_i}, d_j)|$$
$$= | d(c_i, c_j^{\min A_j}) - d(d_i^{\min A_j}, d_j) + d(d_i^{\min A_j}, d_j) - d(y_i, d_j) + d(y_i, d_j) - d(d_i^{\min A_i}, d_j)|$$
$$\leq | d(c_i, c_j^{\min A_j}) - d(d_i^{\min A_j}, d_j)| + |d(d_i^{\min A_j}, d_j) - d(y_i, d_j)| + |d(y_i, d_j) - d(d_i^{\min A_i}, d_j)|$$
$$\leq | d(c_i, c_j^{\min A_j}) - d(d_i^{\min A_j}, d_j)| + d(d_i^{\min A_j}, y_i) + d(y_i, d_i^{\min A_i})$$
The last two terms are less than $g(i)$ since $A_j \subseteq A_i$ and (iii). The first term is less than $g(i) + g(j)$ since these come from Player 2 winning response in the appropriate play of $G^{g,\bar{a}_{\min A_j}, \bar{b}_{\min A_j}}_\alpha$. 
$$ \leq g(i) + g(j) + g(i) + g(i) = 3 g(i) + g(j) = 3 \frac{f(i)}{3} + \frac{f(j)}{3} < f(i) + f(j)$$ 
So the above describes a winning strategy for $G^{f,\bar{a},\bar{b}}_\alpha$. $\bar{a} \sim_\alpha^f \bar{b}$. This completes the proof.
\end{proof}

These are all the results in pure metric space theory that will be needed.

\section{Admissibility}\label{admissibility}
In order to establish bounds on the Scott rank that come from recursion theory or constructibility theory, one needs to look at admissible sets.

$\KP$ is an axiom system in the language $\{\dot \in\}$ where $\dot \in$ is a binary relation symbols. $\KP$ is a weak axiom system for set theory: It includes the basic axioms of set theory such as pairing, union, foundation, and others. The more distinguishing axioms schemes are $\Delta_1$-separation and $\Sigma_1$-replacement. An admissible set is a transitive set $A$ so that $(A,\in) \models \KP$. See \cite{Admissible-Sets-and-Structures}, \cite{Jensen-Model-Existence-Theorem}, \cite{Admissible-Sets}, or \cite{Model-Theory-for-Infinitary-Logic} for more on admissible sets.

As usual in set theory, for emphasis, $V$ will refer to the real universe.

\Begin{definition}{admissible ordinals}
An ordinal $\alpha$ is an admissible ordinal if and only if there is an admissible set $\CA$ so that $A \cap \text{ON} = \alpha$. If $x \in \bairespace$, then $\alpha$ is an $x$-admissible ordinal if and only if there is an admissible set $\CA$ with $x \in A$ so that $A \cap \mathrm{ON} = \alpha$.

For any $x \in \bairespace$, $\omega_1^x$ is the smallest ordinal $\alpha$ so that $L_\alpha(x) \models \KP$.
\end{definition}

\Begin{fact}{basic admissibility results}
An ordinal $\alpha$ is an $x$-admissible ordinal if and only if $L_\alpha(x) \models \mathsf{KP}$. 

If $x \in \bairespace$, $L_{\omega_1^x}(x)$ is the smallest admissible set containing $x$. The reals of $L_{\omega_1^x}(x)$ are the $x$-hyperarithmetic elements. 

$\omega_1^x$ is the supremum of the $x$-hyperarithmetic ordinal as well as the supremum of the $x$-recursive ordinals.
\end{fact}

An important fact about $\KP$ is that the well-founded part of any model of $\KP$ is a model of $\KP$:

\Begin{fact}{truncation lemma}
(Truncation Lemma) Let $\CB = (B,\dot \in^\CB) \models \mathsf{KP}$. Let $\mathrm{WF}(\CB)$ be the collection of $\dot \in^\CB$ well-founded (in $V$) elements of $B$. Then $\mathrm{WF}(\CB) \models \KP$. Hence the Mostowski collapse of $\mathrm{WF}(\CB)$ is an admissible set. 
\end{fact}

\begin{proof}
See \cite{Admissible-Sets-and-Structures}, Lemma II.8.4.
\end{proof}

Let $\SCRL$ be any language. Let $\SCRL_{\infty,\omega}$ denote infinitary logic in the language $\SCRL$. If $\CA$ is an admissible set, then $\SCRL_\CA = (\SCRL_{\infty,\omega})^\CA$. This is the admissible fragment of $\SCRL_{\infty,\omega}$ determined by $\CA$. $\SCRL_\CA$ is a countable admissible fragment if $\CA$ is a countable admissible set. (See \cite{Admissible-Sets-and-Structures} or \cite{Model-Theory-for-Infinitary-Logic} for more information.) The following will be a useful method of constructing admissible sets:

\Begin{fact}{admissible model existence theorem}
(Jensen's model existence theorem) Let $\CA$ be a countable admissible set. Let $\SCRL$ be a language which is $\Delta_1$ definable over $\CA$ and contains a binary relation symbol $\dot \in$ and constant symbols $\hat{a}$ for each $a \in A$. Let $T$ be a consistent theory in the countable admissible fragment $\SCRL_\CA$ which is $\Sigma_1$ definable over $\CA$ and contains the following sentences:

(I) $\mathsf{KP}$

(II) For each $a \in A$, $(\forall v)(v \dot \in \hat{a} \Leftrightarrow \bigvee_{z \in a} v = \hat{z})$. 

\noindent Then there is a $\CB \models T$ so that $\mathrm{WF}(\CB)$ is transitive, $\CA \subseteq \CB$, and $\mathrm{ON} \cap B = \mathrm{ON} \cap A$. 
\end{fact}

\begin{proof}
See \cite{Admissible-Sets} Section 4, Lemma 11 and \cite{Jensen-Model-Existence-Theorem}. Arguments using some form of this fact appear in the proof of Sacks theorem about countable admissible ordinals by Friedman and Jensen. A similar fact is used in Grilliot's omitting type proof of this theorem of Sacks (see \cite{Model-Theory-for-Infinitary-Logic}, Theorem 15).
\end{proof}

\Begin{fact}{barwise compactness}
(Barwise Compactness) Let $\CA$ be a countable admissible set. Let $\SCRL$ be a language which is $\Delta_1$ over $\CA$. Let $T$ be a theory in the countable admissible fragment $\SCRL_\CA$ which is $\Sigma_1$ over $\CA$. If every $F \subseteq T$ so that $F \in A$ is consistent, then $T$ is consistent. 
\end{fact}

\begin{proof}
See \cite{Admissible-Sets-and-Structures} Theorem III.5.6, \cite{Admissible-Sets}, Section 4, Corollary 8, or \cite{Jensen-Model-Existence-Theorem}.
\end{proof}

Now returning back to metric spaces. The following fact expresses how to use an ill-founded ordinal $\alpha$ to play the game $G_\alpha^{f,\bar{a},\bar{b}}$ forever in $V$ to produce an $\SCRU$-automorphism:

\Begin{fact}{ill founded model and automorphisms}
Let $\bar{a}$ and $\bar{b}$ be tuples in $\mathcal{C}(\CM)$ (that is, tuples of $\CM$-Cauchy sequences) of the same length. If there exists an ill-founded model $\CA$ of $\mathsf{ZFC - P}$ with $\mathrm{WF}(\CA)$ transitive so that $\CM,\bar{a},\bar{b} \in A$, $\CA \models \bar{a} \sim_\alpha^f \bar{b}$ for some $f \in \mathrm{REC}$, and $\alpha$ is an ordinal of $\CA$ which is ill-founded (in $V$), then there is a $\SCRU$-automorphism of $\CC(\CM)$ taking $\bar{a}$ to $\bar{b}$. 
\end{fact}

\begin{proof}
Note that since $\CM$ is a metric space on $\omega$, the sets $\CM$, $\bar{a}$, and $\bar{b}$ belong to $\mathrm{WF}(\CA)$. In $\CA$, fix a winning strategy for $G^{f,\bar{a},\bar{b}}_\alpha$ for Player 2.

Let $\Phi : \omega \rightarrow \omega$ be a surjection such that for all $k \in \omega$, $\Phi^{-1}(\{k\})$ is infinite. Let $c_{2i} = \Phi(i)$. Let $d_{2i + 1} = \Phi(i)$. (Recall that $\CM$ is assumed to be a metric space with domain $\omega$.)

Since $\alpha$ is ill-founded, externally in $V$, choose in $V$ an infinite $\CA$-decreasing sequence of $\CA$-ordinals $(\alpha_n)_{n \in \omega}$: that is, for all $n \in \omega$, $\CA \models \alpha_{n + 1} < \alpha_n$. 

Using the winning strategy for Player 2, play as follows:
$$\begin{array}{c | c c c c c c c c c c c c c }
\bar{a} & (\alpha_0,c_0) & {} & {} & c_1 & (\alpha_2,c_2) & {} & {} & c_3 & \dots & (\alpha_{k - 1},c_{k - 1}) & {}  \\
\hline
\bar{b} & {} & d_0 & (\alpha_1, d_1) & {} & {} & d_2 & (\alpha_3, d_3) & {} & \dots  & {} & d_{k - 1}
\end{array}
$$
where for even $i$, $c_i$ are defined above  and for odd $i$, $d_i$ are defined above. For even $i$, $d_i$ comes from the response of Player 2. Similarly for odd $i$, $c_i$ comes from the response of Player 2. Since distance in $\CC(\CM)$ can be expressed as a $\Delta_1$ statement in $\KP$, $\Delta_1$ absoluteness from $\CA$ down to $\mathrm{WF}(\CA)$ and then up into $V$ shows that distance is computed correctly in $\CA$. Hence in $V$, Player 2 has not lost any finite play of $G^{f,\bar{a},\bar{b}}_\alpha$ described above. Since $(\alpha_n)_{n \in \omega}$ is infinite decreasing, the game can always be extended. Playing the game forever in $V$ produces a sequence $(c_n)_{n \in \omega}$ and $(d_n)_{n \in \omega}$ so that each finite portion of the sequence fits into the above play where Player 2 has not lost. 

Let $\Psi : \omega \rightarrow \omega$ be defined by $\Psi(k) = d_k$. Let $\Lambda(k) = c_k$. 

Now to define a map $\Xi : \CC(\CM) \rightarrow \CC(\CM)$: Let $e \in \CC(\CM)$. Let $e = (e^n)_{n \in \omega}$ be some $\CM$-Cauchy sequence representing $e$. Let $\ell : \omega \rightarrow \omega$ be a strictly increasing sequence so that for all $n$, $\Lambda(\ell(n)) = e^n$. Let $\Xi(e)$ be the element of $\CC(\CM)$ represented by the $\CM$-Cauchy sequence $(\Psi(\ell(n)))_{n \in \omega}$.

It straightforward (using argument similar to those of Section \ref{bounds for compact polish metric spaces}) to check that $\Xi$ is well-defined, that is, it does not depend on the Cauchy representation of $e$ or the choice of $\ell$. Using the definition of $G^{f,\bar{a},\bar{b}}_\alpha$, one can check that $\Xi$ is a $\SCRU$-homomorphism and that $\bar{a}$ is mapped to $\bar{b}$. By how $\Phi$ was chosen, one can show that $\Xi$ is actually surjective. Hence $\Xi$ is a $\SCRU$-automorphism taking $\bar{a}$ to $\bar{b}$. (It should be noted that the fact that $\Phi^{-1}(\{i\})$ is infinite for each $i \in \omega$ is important for establishing these properties.)
\end{proof}

In the following, $\bar{a}$ and $\bar{b}$ are considered as tuples of $\CM$-Cauchy sequences. Since $\CM$ is a metric space on $\omega$, $\bar{a}$ and $\bar{b}$ may be coded as elements of $\bairespace$. 

\Begin{fact}{bounds on metric rank of tuple}
If there is no $\SCRU$-automorphism taking $\bar{a}$ to $\bar{b}$, then $\RR(\bar{a},\bar{b}) < \omega_1^{\CM \oplus \bar{a} \oplus \bar{b}}$ and in particular, $\SR(\bar{a},\bar{b}) < \omega_1^{\CM \oplus \bar{a} \oplus \bar{b}}$. 

(\cite{Scott-Rank-of-Polish-Metric-Spaces} Doucha) If there is no $\SCRU$-automorphism taking $\bar{a}$ to $\bar{b}$, then $\SR(\bar{a},\bar{b})$ is countable. Therefore, the Scott rank of a Polish metric space is at most $\omega_1$.
\end{fact}

\begin{proof}
Suppose $\RR(\bar{a},\bar{b}) \geq \omega_1^{\CM \oplus \bar{a} \oplus \bar{b}}$. 

Let $\CA = L_{\omega_1^{\CM \oplus \bar{a} \oplus \bar{b}}}(\CM \oplus \bar{a} \oplus \bar{b})$. $\CA$ is a countable admissible set. 

Let $\SCRL$ be a language consisting of the following:

(i) A binary relation symbol $\dot \in$. 

(ii) For each $a \in A$, a constant symbol $\hat{a}$. 

\noindent $\SCRL$ is a language which is $\Delta_1$-definable in $\CA$. 

Now let $T$ be the theory in the countable admissible fragment $\SCRL_\CA$ consisting of sentences indicated below:

(I) $\mathsf{ZFC - P}$. 

(II) For each $a \in A$, ``$(\forall v)(v \dot \in \hat{a} \Leftrightarrow \bigvee_{z \in a} v = \hat{z})$''. 

(III) For each $\alpha < \omega_1^{\CM \oplus \bar{a} \oplus \bar{b}}$, ``$\bigwedge_{f \in \mathrm{REC}} \hat{\bar{a}} \sim_\alpha^f \hat{\bar{b}}$''. 

\noindent $T$ is $\Sigma_1$ definable in $\CA$. 

$T$ is consistent. To see this, consider the structure $\CB$ defined by: Let its domain be $B = H_{\aleph_1}$, the collection hereditarily countable sets. Let $\dot \in = \in \upharpoonright H_{\aleph_1}$. For each $a \in A$, let $\hat{a}^{\CB} = a$. $\CB \models T$ since it was assumed that $\RR(\bar{a},\bar{b}) \geq \omega_1^{\CM \oplus \bar{a} \oplus \bar{b}}$. 

By Fact \ref{admissible model existence theorem}, there is a $\SCRL$-structure $\CB$ so that $\CB \models T$, $\text{WF}(\CB)$ is transitive, $\text{ON} \cap B = \text{ON} \cap A = \omega_1^{\CM \oplus \bar{a} \oplus \bar{b}}$, and $\CA$ is an end extension of $\CB$. $\CB$ must be ill-founded since no transitive set of ordinal height $\omega_1^{\CM\oplus\bar{a}\oplus\bar{b}}$ containing $\CM\oplus\bar{a}\oplus\bar{b}$ can be a model of $\mathsf{ZFC - P}$. Hence by (III), there must be some illfounded $\CB$-ordinal $\beta$ so that for any $f \in \text{REC}$, $\CB \models \bar{a} \sim_\beta^f \bar{b}$.

Fact \ref{ill founded model and automorphisms} shows that there is a $\SCRU$-automorphism taking $\bar{a}$ to $\bar{b}$. Contradiction.
\end{proof}

Next, it will be shown that if $\RR(\bar{a},\bar{b}) \geq \omega_1^\CM$, then $\RR(\bar{a},\bar{b})$ is a limit of elements $(\bar{e},\bar{f})$ so that there is a $\SCRU$-automorphism of $\CC(\CM)$ taking $\bar{e}$ to $\bar{f}$.

In fact, in \cite{Scott-Rank-of-Polish-Metric-Spaces} Proposition 2.4, it is shown that certain points have the property that every open neighborhood contains a perfect set of $(\bar{e},\bar{f})$ so that there is a $\SCRU$-automorphism of $\CC(\CM)$ taking $\bar{e}$ to $\bar{f}$. It can also be shown that if $\RR(\bar{a},\bar{b}) \geq \omega_1^\CM$, then every open set containing $(\bar{a},\bar{b})$ has a perfect set of such $(\bar{e},\bar{f})$. However, this fact is not necessary for producing the effective bound on Scott rank. The reader may choose to skip all the comments about perfect sets in the following two results.

If one wants the perfect set result, one will need the following effective perfect set theorem:

\Begin{fact}{harrison effective perfect set}
(Harrison) Let $r \in \bairespace$. Suppose $X$ is a $r$-recursively presented Polish space. If a $\Sigma_1^1(r)$ set $A$ contains a member which is not $\Delta_1^1(r)$ (i.e. $r$-hyperarithmetic), then $A$ contains a perfect subset.
\end{fact}

\begin{proof}
See \cite{Descriptive-Set-Theory} for more on recursively presented Polish space and \cite{Descriptive-Set-Theory} Theorem 4F.1.
\end{proof}

\Begin{fact}{high scott rank limit of automorphism moving points}
Let $\bar{a}$ and $\bar{b}$ be tuples of elements of $\CC(\CM)$ of the same length $p$. Suppose $\RR(\bar{a},\bar{b}) \geq \omega_1^\CM$ and there are no $\SCRU$-automorphisms of $\CC(\CM)$ taking $\bar{a}$ to $\bar{b}$. Then for any $\bar{n} \in {}^{2p}M$ and $m \in \omega$, if $(\bar{a},\bar{b}) \in B_\frac{1}{m}(\bar{n})$ (the open ball around $\bar{n}$ of size $\frac{1}{m}$ in the metric on ${}^{2p}\CC(\CM)$), then there is some $(\bar{e},\bar{f}) \in B_\frac{1}{m}(\bar{n})$ for which there is a $\SCRU$-automorphism taking $\bar{e}$ to $\bar{f}$.

In fact, there is a perfect set of such $(\bar{e},\bar{f})$. 
\end{fact}

\begin{proof}
Fix $m \in \omega$ and $\bar{n} \in {}^{2p}M$ so that $(\bar{a},\bar{b}) \in B_{\frac{1}{m}}(\bar{n})$. 

Let $\CA = L_{\omega_1^\CM}(\CM)$. $\CA$ is a countable admissible set. 

(Some remarks before continuing: Since $\bar{a}$ and $\bar{b}$ are tuples of $\CM$-Cauchy sequences, they are coded by reals. As it will be shown below, $\bar{a}$ and $\bar{b}$ can not belong to $\CA$. Thus one cannot mention $\bar{a}$ or $\bar{b}$ in any countable fragment associated to the admissible set $\CA$. However, $\bar{n}$ is a tuple of elements of $\CM$ which is essentially a tuple of integers (since $\CM$ was assumed to be a metric space on $\omega$). Thus $\bar{n}$ belongs to any admissible set. One is permitted to refer to $\bar{n}$. Although $\bar{a}$ and $\bar{b}$ cannot be mentioned in the theory, these elements will be used to (externally in $V$) verify the consistency of the theory. The details follow as the proof resumes below.)

Let $\SCRL$ be a language consisting of the following:

(i) A binary relation symbol $\dot \in$. 

(ii) For each $a \in A$, a constant symbol $\hat{a}$. 

(iii) Two new constant symbols $\dot{\bar{e}}$ and $\dot{\bar{f}}$. 

\noindent $\SCRL$ is a language which is $\Delta_1$ definable in $\CA$. 

Now let $T$ be the theory in the countable admissible fragment $\SCRL_\CA$ consisting of the sentences indicated below:

(I) $\mathsf{ZFC - P}$.

(II) For each $a \in A$, ``$(\forall v)(v \dot\in \hat{a} \Leftrightarrow \bigvee_{z \in a} v = \hat{z})$''. 

(III) ``$d_{2p}(\hat{\bar{n}}, (\dot{\bar{e}}, \dot{\bar{f}})) < \frac{1}{m}$'' where $d_{2p}$ is the metric on ${}^{2p}\CC(\CM)$. 

(IV) For each $\alpha < \omega_1^{\CM}$, ``$\bigwedge_{f \in \mathrm{REC}} \dot{\bar{e}} \sim_\alpha^f \dot{\bar{e}}$''. 

\noindent If one want the perfect set version of this result, add on the following 

(V) For all $\alpha< \omega_1^{\CM}$, ``$(\dot{\bar{e}},\dot{\bar{f}}) \notin L_{\alpha}(\hat{\CM})$''. 

\noindent In either case, $T$ is $\Sigma_1$ definable in $\CA$. 

$T$ is consistent. To see this: Consider the following structure $\CB$. Its domain is $B = H_{\aleph_1}$. $\dot \in^{\CB} = \in \upharpoonright H_{\aleph_1}$. For each $a \in A$, let $\hat{a}^\CB = a$. Let $\dot{\bar{e}}^\CB = \bar{a}$ and $\dot{\bar{f}}^\CB = \bar{b}$. $\CB \models T$ since $\RR(\bar{a},\bar{b}) \geq \omega_1^\CM$. 

For those interested in the perfect set version, to see (V), note that $\omega_1^{\CM \oplus \bar{a}\oplus\bar{b}} > \omega_1^\CM$. If not, then $\RR(\bar{a},\bar{b}) \geq \omega_1^{\CM} \geq \omega_1^{\CM \oplus \bar{a} \oplus \bar{b}}$. By Fact \ref{bounds on metric rank of tuple}, there is a $\SCRU$-automorphism taking $\bar{a}$ to $\bar{b}$. This contradicts the assumption on $(\bar{a},\bar{b})$. Since $\omega_1^{\CM \oplus \bar{a} \oplus \bar{b}} > \omega_1^{\CM}$, $\CM \oplus \bar{a} \oplus \bar{b}$ can not belong to any admissible set of ordinal height $\omega_1^\CM$. In particular, $(\bar{a},\bar{b}) \notin L_{\omega_1^{\CM}}(\CM)$. This shows the model $\CB$ satisfies (V). 

By Fact \ref{admissible model existence theorem}, there exists some model $\CB \models T$ so that $\mathrm{WF}(\CB)$ is transitive and $\mathrm{ON} \cap B = \mathrm{ON} \cap A$. Let $\bar{e} = \dot{\bar{e}}^{\CM}$ and $\bar{f} = \dot{\bar{e}}^{\CM}$. As before, $\CB$ must be ill-founded. In $\CB$, there is some ill-founded $\CB$-ordinal $\alpha$ so that $\CB \models \bar{e} \sim_\alpha^f \bar{f}$, for any $f \in \mathrm{REC}$. So by Fact \ref{ill founded model and automorphisms} in $V$, there is a $\SCRU$-automorphism taking $\bar{e}$ to $\bar{f}$. By absoluteness, $(\bar{e},\bar{f}) \in B_\frac{1}{m}(\bar{n})$. Hence $(\bar{e},\bar{f})$ is the desired element. 

For the perfect set version, note that the set of $C$ of elements $(\bar{u},\bar{v}) \in B_\frac{1}{m}(\bar{n})$ so that there is a $\SCRU$-automorphism taking $\bar{u}$ to $\bar{v}$ is a $\Sigma_1^1(M)$ set. By (V), the element $(\bar{e},\bar{f})$ produced above is not in $L_{\omega_1^{\CM}}(\CM)$, so in particular not $\lborel(\CM)$. Hence $C$ must contain a perfect subset by Fact \ref{harrison effective perfect set}.
\end{proof}

\section{Main Results}\label{main results}

\Begin{definition}{rigid structure}
Let $\SCRL$ be a language and let $\CN$ be an $\SCRL$-structure. $\CN$ is \text{rigid} if and only if there are no nontrivial $\SCRL$-automorphisms of $\CN$. 
\end{definition}

\Begin{theorem}{rigid structure low scott rank}
Let $\CM$ be a metric space on $\omega$. Suppose $\CC(\CM)$ is a rigid Polish metric space. Then $\SR(\CC(\CM)) < \omega_1^\CM$. 
\end{theorem}

\begin{proof}
Suppose $\SR(\CC(\CM)) \geq \omega_1^\CM$. This means for each $\alpha < \omega_1^\CM$, there is some $\bar{a}_\alpha$ and $\bar{b}_\alpha$ so that $\bar{a}_\alpha \neq \bar{b}_\alpha$ (as elements of $\CC(\CM)$) and $\bar{a} \sim_\alpha^f \bar{b}$ for any $f \in \mathrm{REC}$. Note that for $\alpha \neq \beta$, the length of $\bar{a}_\beta$ and $\bar{a}_\beta$ may not be the same.

Let $\CA$ be $L_{\omega_1^\CM}(\CM)$. Let $\SCRL$ be a language consisting of:

(i) A binary relation symbol $\dot \in$. 

(ii) For each $a \in A$, a constant symbol $\hat{a}$. 

(iii) Three new constant symbols $\dot n$, $\dot e$, and $\dot f$. 

\noindent $\SCRL$ is $\Delta_1$ definable over $\CA$. 

If $n \in \omega$ and $r \in \bairespace$, let $c_n(r)$ denote the element of ${}^n(\bairespace)$ coded by $r$ (under some fixed recursive coding of $n$-tuples of reals by a single real).

Let $T$ be the following theory in the countable admissible fragment $\SCRL_\CA$ consisting of the sentences indicated below:

(I) $\mathsf{ZFC - P}$

(II) For each $a \in A$, ``$(\forall v)(v \dot \in \hat{a} \Leftrightarrow \bigvee_{z \in a} v = \hat{z})$''. 

(III) ``$\dot n \dot \in \hat{\omega}$''. ``$\dot e$ and $\dot f$ are functions from $\hat{\omega}$ to $\hat{\omega}$''. 

(IV) ``$c_{\dot n}(\hat{e})$ and $c_{\dot n}(\hat{f})$ are tuples of $\CM$-Cauchy sequences''. ``$c_{\dot n}(\hat{e}) \neq c_{\dot n}(\hat{f})$ as $\CM$-Cauchy sequences''.

(V) For each $\alpha < \omega_1^\CM$, ``$\bigwedge_{f \in \mathrm{REC}} c_{\dot n}(\dot e) \sim_\alpha^f c_{\dot n}(\dot f)$''. 

\noindent $T$ is $\Sigma_1$ definable over $\CA$. 

To see that $T$ is consistent, one needs to use Barwise compactness. Let $F \subset T$ be such that $F \in A$. Since $F \in \CA$, there is an $\alpha < \omega_1^\CM$ that bounds all the $\beta$'s that appear in statements of type (IV). Consider the model $\CB$ defined by: Its domain is $B = H_{\aleph_1}$. $\dot \in^\CB = \in \upharpoonright \CB$. For each $a \in A$, $\hat{a}^\CB = a$. Let $\dot n^\CB = |\bar{a}_\alpha|$. Let $e$ and $f$ be two reals so that $c_{|\bar{a}_\alpha|}(e) = \bar{a}_\alpha$ and $c_{|\bar{a}_{\alpha}|}(f) = \bar{b}_\alpha$. It is clear that $\CB \models F$. By Barwise compactness, $T$ is consistent. 

Now Fact \ref{admissible model existence theorem} gives a model $\CB \models T$ so that $\mathrm{WF}(\CB)$ is transitive, $\CA \subseteq \CB$, and $\mathrm{ON} \cap \CB = \mathrm{ON} \cap \CA = \omega_1^\CM$. Let $\bar{e} = c_{\dot n^\CB}(\dot e^\CB)$, $\bar{f} = c_{\dot n^\CB}(\dot f^\CB)$. Since $\CM$ is a metric space on $\omega$, all $\CM$-Cauchy sequences in $\CB$ belong to $\mathrm{WF}(\CB)$. Hence $\bar{e},\bar{f} \in \mathrm{WF}(\CB)$. By $\Delta_1$-absoluteness from $\CB$ down to $\mathrm{WF}(\CB)$ and then up to $V$, one can show that for all $\alpha < \omega_1^\CM$, $\bar{e} \sim_\alpha^f \bar{f}$ for all $f \in \mathrm{REC}$. Also by $\Delta_1$-absoluteness, $\bar{e} \neq \bar{f}$. However $\CM \oplus\bar{e} \oplus \bar{f}$ is in the admissible set $\mathrm{WF}(\CB)$ (by Fact \ref{truncation lemma}) which has ordinal height $\omega_1^\CM$. Hence $\omega_1^{\CM \oplus \bar{e} \oplus \bar{f}} = \omega_1^\CM$. Fact \ref{bounds on metric rank of tuple} implies there is a $\SCRU$-automorphism taking $\bar{e}$ to $\bar{f}$. This contradicts the assumption that $\CC(\CM)$ is a rigid metric space.
\end{proof}

\Begin{theorem}{proper polish space bound on scott rank}
Let $\CM$ be a metric space on $\omega$. Suppose $\CC(\CM)$ is a proper Polish metric space. Then $\SR(\CM) \leq \omega_1^{\CM} + 1$. 
\end{theorem}

\begin{proof}
Suppose not, then there exists some tuples $\bar{a}$ and $\bar{b}$ of elements of $\CC(\CM)$ of the same length so that $R(\bar{a},\bar{b}) \geq \omega_1^\CM$ but there is no $\SCRU$-automorphism taking $\bar{a}$ to $\bar{b}$. By Fact \ref{high scott rank limit of automorphism moving points}, there exist a sequence $(\bar{a}_n,\bar{b}_n)_{n \in \omega}$ so that $(\bar{a},\bar{b})$ is its limit and for all $n \in \omega$, there is a $\SCRU$-automorphism taking $\bar{a}_n$ to $\bar{b}_n$. Let $\alpha$ be an ordinal greater than $\omega_1^{\CM \oplus \bar{a} \oplus \bar{b}}$. The existence of these automorphisms implies that $\bar{a}_n \sim^f_\alpha \bar{b}_n$ for all $n \in \omega$ and any $f \in \mathrm{REC}$. Then Fact \ref{proper polish space limit certain rank implies certain rank} implies that $\bar{a} \sim_\alpha^f \bar{b}$. However since $\alpha > \omega_1^{\CM \oplus \bar{a} \oplus \bar{b}}$, Fact \ref{bounds on metric rank of tuple} implies that there is a $\SCRU$-automorphism taking $\bar{a}$ to $\bar{b}$. Contradiction.
\end{proof}

\bibliographystyle{amsplain}
\bibliography{references}

\end{document}